\documentclass[12pt,a4paper]{article}
\pdfoutput=1

\usepackage{amsbsy,amsfonts,amsmath}
\usepackage[dvips]{graphicx}
\usepackage{newlfont}
\usepackage{color}
\usepackage{url}
\usepackage{caption}
\usepackage{subcaption}

\usepackage{pgfplots}
\usepackage{tikz}
\usepackage{setspace}
\pgfplotsset{compat=1.12}
\usepackage{lscape}
\usepackage{booktabs,multicol, multirow}

\newcommand{\subjto}{\ensuremath{\text{subject to}}}

\definecolor{orange}{rgb}{1,0.5,0}

\def\keywords{\vspace{-.3em}
 \if@twocolumn
  \small\it Keywords\/\bf---$\!$%
 \else
\begin{center}\small\bf
Keywords\end{center}\quotation\small  \fi}
\def\endkeywords{\vspace{0.6em}\par\if@twocolumn\else\endquotation\fi
 \normalsize\rm}

\title{Metro stations as crowd-shipping catalysts: \\ an empirical and computational study}
\author{Carlo Filippi\thanks{Department of Economics and Management, University of Brescia, Contrada S. Chiara 50, 25122 Brescia, Italy - email:\{carlo.filippi; francesca.plebani\}@unibs.it} \quad Francesca Plebani\footnotemark[1]}

\date{}
\begin{document}

\bibliographystyle{abbrv}

\maketitle

\begin{abstract}
Crowd-shipping is a promising shared mobility service that involves the delivery of goods using non-professional shippers. This service is mainly intended to reduce congestion and pollution in city centers but, as some authors observe, in most crowd-shipping initiatives the crowd rely on private motorized vehicles and hence the environmental benefits could be small, if not negative. Conversely, a crowd-shipping service relying on public transport should maximize the environmental benefits. Motivated by this observation, in this study we assess the potentials of crowd-shipping based on metro commuters in the city of Brescia, Italy. Our contribution is twofold. First, we analyze the results of a survey conducted among metro users to assess their willingness to act as crowd-shippers. The main result is that most young commuters and retirees are willing to be crowd-shippers even for a null reward. Second, we assess the potential economic impact of using metro-based crowd-shipping coupled with a traditional home delivery service. To this end, we formulate a variant of the VRP model where the customers closest to the metro stations may be served either by a conventional vehicle or by a crowd-shipper. The model is implemented using Python with Gurobi solver. A computational study based on the Brescia case is performed to get insights on the economic advantages that a metro-based crowd delivery option may have for a retailing company.
\end{abstract}

\noindent

\section{Introduction}
\label{sec:intro}
E-commerce, home-delivery, same-day delivery continue to exhibit a remarkable growth year after year. These phenomena are symptoms of a profound change in retailing and in consumer behavior and involve huge challenges for last-mile logistics management. 
On one side, direct deliveries are likely to cause smaller loads and more frequent shipments, involving lower freight consolidation. Consequently, logistics companies face an increment in marginal transportation costs coupled with an increasing competitiveness, and are forced to experiment with innovative and cost-effective delivery systems. On the other side, direct deliveries generate a worsening of traffic and parking conditions in urban areas, while public authorities are struggling to reduce traffic congestion and pollution in our cities. Thus, the urge of efficient and sustainable delivery modes is experienced by private and public stakeholders.

Crowd-logistics is a concept that may help in this context. Crowd-logistics ``designates the outsourcing of logistics services to a mass of actors, whereby the coordination is supported by a technical infrastructure. The aim of crowd-logistics is to achieve economic benefits for all stake- and shareholders'' \cite{2015-MehFreTeu}. Crowd-logistics encompasses different sharing-economy models, including crowd-shipping, where individuals traveling to a certain region are willing to perform deliveries on their way, and professional shippers may rely on them to carry out part of their shipments. The integration of personal and freight transport is allowed by on-line platforms that, exploiting mobile communication technologies and global positioning systems (GPS), are able to match in real time delivery demand with transportation supply. 

Crowd-shipping may be implemented in different ways. In most cases, crowd-shippers are ``occasional drivers'', that pick up a parcel and deliver it to the final customer by using privately owned means of transportation. Examples of initiatives of this type are ``Spark Delivery'' by Walmart, ``Amazon Flex'', ``Deliv'' in several U.S. cities, ``Trunkrs'' in The Netherlands, and ``Renren Kuaidi'' in China \cite{2019-ArsAgaKroZui}. In this implementation, crowd-shippers are willing to pick up a parcel along a journey that they are already making and deliver it making a limited detour. 

While the use of a motorized private vehicle is not an obstacle to implement a cost-effective shipment service, it could even worsen congestion and pollution in the urban complex. Hence, while the positive effects for companies in terms of cost and image are interesting, the positive effects for the quality of life are questionable \cite{2018-QiLiLiuShe}. A more sustainable implementation of crowd-shipping might rely on public transport, i.e., on the trips of commuters using subways, light rails, buses. In this case, any parcel shifted from professional delivery to crowd-shipment should have a positive effect on the environment. Real crowd-shipment initiatives based of public service are harder to find. A remarkable example is the ``Dabbawala'' lunchbox delivery and return system in India, especially Mumbai, which is based on railway trains and bycicles \cite{2013-BaiMac}.

In this paper, we analyze the potentials of a crowd-shipping initiative based on public transport for a European medium-sized city.
More specifically, we consider the case of Brescia, a city of around 200,000 inhabitants situated in northern Italy. Brescia is County Seat of the largest province in Lombardy region and an academic center. From a few years, Brescia owns a fully automated metro system with a single L-shaped line, connecting the south-east part of the city to the northern part. On the Brescia case, we formulate two research questions:
\begin{itemize}
\item[Q1] Are metro commuters willing to operate as crowd-shippers, adding a pickup and delivery routine to their usual path, possibly after a small reward?
\item[Q2] Is it economically interesting for a professional shipper the chance of using metro commuters as occasional delivery persons, possibly paying them a small reward?  
\end{itemize}
To respond to the first question, we handed out a questionnaire to a sample of Brescia metro users and we analyzed the obtained answers. It turns out that metro users are very sensitive to environmental sustainability, and that most commuters are willing to be crowd-shippers, even for a null reward, especially if they are young or aged.
To respond to the second question, we formulated and implemented an optimization model to assess the benefits for a private professional shipper to use crowd-shipping based on metro commuters. More specifically, we consider a vehicle routing problem (VRP) in which the delivery locations reasonably close to metro stations can be served either by a professional shipper or by a crowd-shipper (metro user), possibly paying a fixed reward. 
We use the model to run an extensive simulation involving Brescia municipality and its metro system. We obtain managerial insights of the profitability of a crowd-shipping initiative based on public transport.

The paper is organized as follows. In Section~\ref{sec:litrev}, we analyze the relevant literature. In Section~\ref{sec:survey}, we describe the survey conducted in Brescia and discuss its results. In Section~\ref{sec:prodes}, we design an optimization model for the simulation of the usage of metro-based crowd-shipping in Brescia, we develop the case study and analyze it using the proposed model. Finally, Section~\ref{sec:concl} is devoted to some conclusions.

\section{Literature review}
\label{sec:litrev}

A recent and thorough survey on crowd-shipping is proposed by Le et al. \cite{2019-LeStaVanUkk}, who review current practice, academic research, and empirical case studies according to three ``pillars'': supply, demand, and operations-and-management. 
Crowd-shipping is a special form of \emph{crowdphysics}, a term coined by Sadilek, Krumm, and Horvitz \cite{2013-SadKruHor} to denote ``crowdsourcing tasks that require people to collaborate and synchronize both in time and physical space.'' 
After this seminal contribution, several works dealt with crowd-shipping under different perspectives.

A large part of the literature focuses on the motivations and economical perspectives of crowd-shipping. The goal is usually to help policy-makers to adapt laws and regulations to the sharing economy developments and to provide insights to businesses in which crowd-shipping is compatible with the corporate social responsibility strategy. 
Buldeo Rai et al. \cite{2017-BulVerMerMac} formulate a definition of the wider concept of \emph{crowd-logistics} and identify which factors determine the sustainability potential of crowd-logistics. Their methodology includes literature analysis and interviews with logistics practitioners. 
Punel and Stathopoulos \cite{2017-PunSta} analyze  the determinants of crowd-shipping acceptance among senders, using a discrete choice model. They characterize different preference patterns related to different distance classes for the crowd-shipping service.  
The difference between crowd-shipping users and non-users is studied by Punel, Ermagun, and Stathopoulos \cite{2018-PunErmSta}. Using proportional $t$-test analysis and a binary logit model they conclude that crowd-shipping is more prevalent among young people, men, and full-time employed individuals, while urban areas are preferential for the development of crowd-shipping. 
Ermagun and Stathopoulos \cite{2018-ErmSta} study the supply determinants of crowd-shipping. They exploit a large data set of crowd-shipping requests across the Unites States to correlate the potential of crowd-shipping initiatives to several social and economic variables, giving insights for transportation planners and crowd-shipping companies. In a following paper \cite{2020-ErmSta}, the same authors use the data set to analyse the performance of existing crowd-shipping services and identify levers to improve efficiency. 

A few papers focus on optimization models and algorithms for crowd-shipping. Archetti, Savelsbergh, and Speranza \cite{2016-ArcSavSpe} propose the Vehicle Routing Problem with Occasional Drivers (VRPOD) to model a situation where deliveries to customer's home is made from a given store and can be done either by conventional vehicles or using an ``occasional driver'', i.e., a customer present at the store and willing to deliver a package to another customer, in exchange of a little reward. Starting from the VRPOD, Macrina et al. \cite{2017-MacDiPGueLag} introduce time windows for the visit of customers and propose and validate two alternative optimization models. In a subsequent work \cite{2020-MacDiPGueLap}, the same authors introduce also transshipment nodes, and design a variable neighborhood search heuristic. Raviv and Tenzer (2018) study a crowd-shipping problem in which the deliveries are carried out by external couriers. For a version of the problem without capacity constraints, they design a stochastic dynamic algorithm and a greedy heuristic.

Most of the studies on crowd-shipping assume that the crowd-shippers are drivers that use their private vehicles to offer a shipping service. A few authors consider the possibility to implement crowd-shipping using public transport. Serafini et al. \cite{2018-SerNigGatMar} and Gatta et al. \cite{2019-GatMarNigSer} investigate the willingness to act as crowd-shipper of metro users in the city of Rome, Italy, where Automated Parcel Lockers (APL) located in metro stations are used for pick-up/drop-off operations. They use stated preferences and discrete choice models to analyze the inclination to crowd-shipping. Starting from the same case study, Gatta et al. \cite{2019-GatMarNigPatSer} evaluate the environmental and economic impact of crowd-shipping in urban areas. Finally, Simoni et al. \cite{2020-SimMarGatCla} develop a hybrid dynamic traffic simulation model to analyze the impact of different operational features, like transport crowd-shipping mode and admissible detour length, on congestion and emissions.



\section{Empirical study: the survey}
\label{sec:survey}
In this section, we investigate the propensity of the inhabitants of the city of Brescia towards the crowd-shipping service combined with the usage of the metro. More specifically, we analyze through a survey the willingness of the inhabitants to act as crowd-shippers (supply side) and to purchase a crowd-shipping service (demand side) to get goods delivered/collected in a last mile B2C e-commerce situation. The idea behind the investigation focuses on the opportunity to involve in a crowd-shipping project the population that normally uses the metro to travel downtown, and deliver goods purchased online to customers located near metro stops. The project could involve both home delivery by crowd-shippers to customers living near metro stations and the pickup/release of goods in automated parcel lockers located at the stations. 
Our survey is inspired by an analogous research done by Gatta et al. \cite{2019-GatMarNigSer} for the city of Rome.

%
%


Brescia is a city of about 200,000 inhabitants distributed on a municipal territory of about 90 km$^2$.
Since 2013, the city is served by a metro, an automatic light rail line that connects the northern districts of the city of Brescia to those of the south-east, passing through the historic center. It employs a fully automatic rapid rail transport system, designed and built by Ansaldo-STS, similar to that already realized for the Copenhagen metro and the contemporary M5 metro in Milan. The service consists of a single line with a total of 17 stations \cite{BresciaMob}. The number of metro users has constantly grown from the inauguration to the Covid-19 pandemics, reaching about 18.7 million in 2019 \cite{BresciaMob}.

Our survey was carried out through a questionnaire created with Google Form and submitted online between the end of July 2020 and October 2020, which allowed to collect a total of 392 responses. 

%

\subsubsection*{The sample}
\begin{figure}
  \centering
  \begin{subfigure}{0.475\textwidth}
    \centering
    \includegraphics[width=0.8\linewidth]{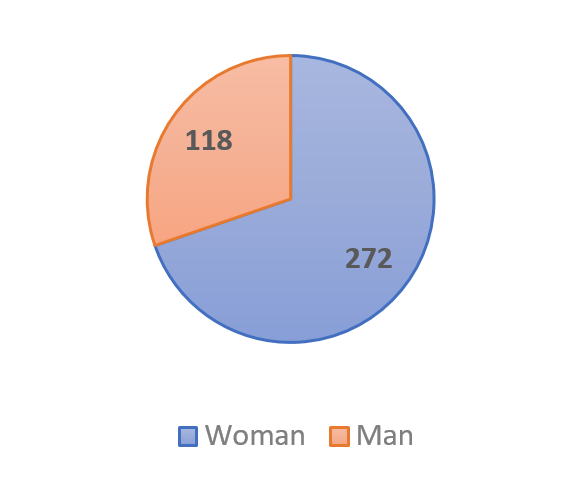}%
		\caption{Gender distribution}
  \end{subfigure}
  \hfill
  \begin{subfigure}{0.475\textwidth}
    \centering
    \includegraphics[width=0.8\linewidth]{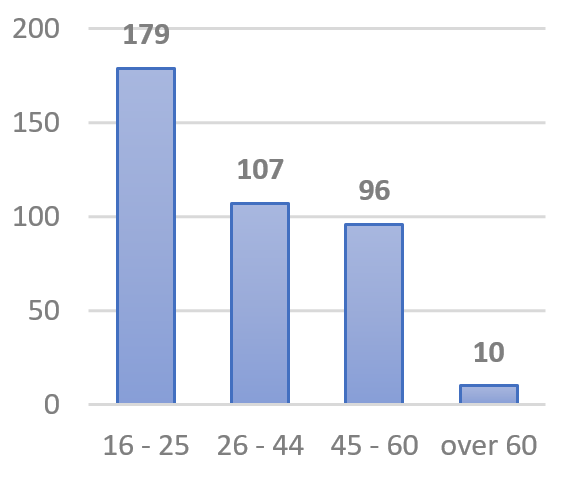}%
		\caption{Age distribution}
  \end{subfigure}\\
  \begin{subfigure}{0.475\textwidth}
    \centering
    \includegraphics[width=0.8\linewidth]{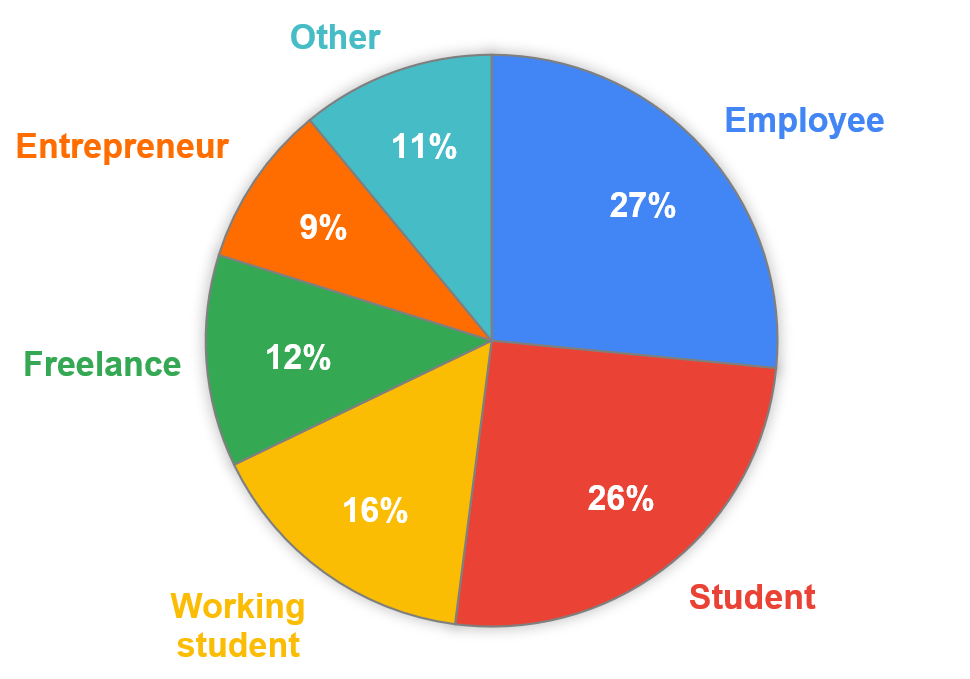}%
		\caption{Professional distribution}
  \end{subfigure}\\	
\caption{Gender, age, and professional distribution of the sample.}
\label{fig:gender-age}
\end{figure}

The first part of the questionnaire was devoted to a simple analysis of the respondents. The obtained distribution by gender, age, and professional status of the sample is depicted in Figure \ref{fig:gender-age}. Notice that the question about gender was not mandatory, but all the participants responded, with 70\% of women and 30\% of men.

\subsubsection*{Interest in the initiative}
\begin{figure}
  \centering
  \begin{subfigure}{0.475\textwidth}
    \centering
    \includegraphics[width=0.8\linewidth]{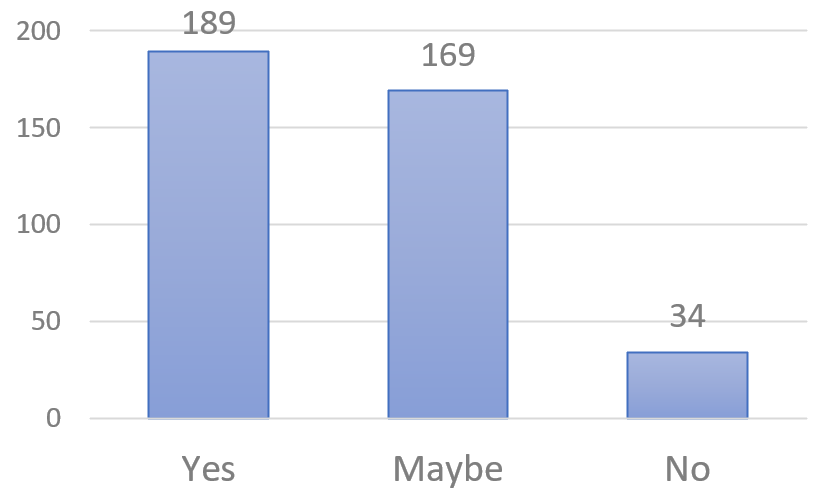}%
		\caption{Future of crowd-shipping in Brescia}
  \end{subfigure}
  \hfill
  \begin{subfigure}{0.475\textwidth}
    \centering
    \includegraphics[width=0.6\linewidth]{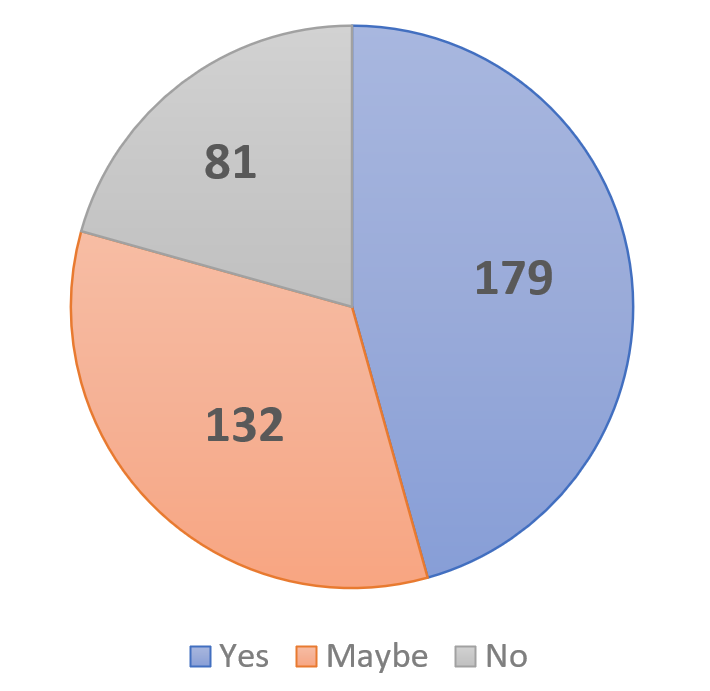}%
		\caption{Willingness to collaborate.}
  \end{subfigure}	
\caption{Future of crowd-shipping and willingness to collaborate.}
\label{fig:future-collaboration}
\end{figure}

The second part of the questionnaire involved the personal interest in a crowd-shipping initiative. In particular, we pose the following questions.
\begin{itemize}
\item \emph{Do you think that crowd-shipping, an initiative that allows for environmental and traffic benefits and a possible decrease in delivery times for goods, could have a future in Italy? (e.g. in big cities like Milan and Rome)} The possible answers was ``yes'', ``maybe'', ``no''. Only 6\% of the interviewees think that crowd-shipping has no future in Italy, against 64\% who are completely optimistic. There is no significant difference between men and women, while age has its influence: the more skeptical are the respondents in age class 26-44 (10\% of ``no'' answers, 53\% of ``yes'' answers); the more optimistic are the respondents over 60, with 9 ``yes'' answers and one ``maybe'' answer.
\item \emph{Do you think crowd-shipping can have a future in Brescia?} Restricting the focus to the city of Brescia, confidence decreases, but only 9\% of respondents think it has no future. The percentage of undecided respondents has risen sharply, to around 43\%, see Figure \ref{fig:future-collaboration}(a). The differences between genders and among age classes are similar to the previous case.
\item \emph{In light of the benefits it could have for Brescia from an environmental point of view, would you take action to see this initiative realized? For example, by talking about it and raising awareness about it among friends and relatives, by participating in initiatives and presentations, by proposing yourself as an occasional shipper}. Here, 21\% of respondents are not willing to collaborate, see Figure \ref{fig:future-collaboration}(b); the answers are partitioned by age class in Table \ref{tab:willingness}. According to the answers, 85\% of the youngest respondents and 90\% of the oldest ones can be considered as potential crowd-shippers. 
This suggests that students and retirees are the most likely to become crowd-shippers.
Finally, we notice that in this case men appear to be more collaborative than women (17\% vs 27\% of ``no'' answers).
\end{itemize}

\begin{table}
  \centering
    \begin{tabular}{p{0.2\textwidth}p{0.1\textwidth}p{0.1\textwidth}p{0.1\textwidth}}
		\toprule
    Age   & Yes & Maybe & No \\
		\midrule
    16 - 25 & 45\%  & 40\%  & 16\% \\
    26 - 44 & 40\%  & 33\%  & 27\% \\
    45 - 60 & 53\%  & 23\%  & 24\% \\
    Over 60 & 50\%  & 40\%  & 10\% \\
		\bottomrule
    \end{tabular}%
  \caption{Willingness to collaborate with crowd-shipping initiatives by age class.}
  \label{tab:willingness}%
\end{table}%

\subsubsection*{Supply-side and demand-side}
In order to understand which and how many of the respondents represent a significant sample for the demand-side and supply-side survey, we have divided the respondents into categories according to the following questions.
\begin{itemize}
\item \emph{Are you a resident of the city of Brescia?} Living in Brescia allows us to identify potential customers of the service with crowd-shipping combined with the use of the metro line. The 104 respondents living in the city of Brescia represent the reference sample for the demand-side survey.
\item \emph{Do you use the Brescia metro regularly during the year?} 288 out of 392 respondents stated that they use the metro only occasionally; therefore, they cannot represent the supply-side reference sample, which is therefore composed of the remaining 104 respondents.
\end{itemize}

The cross-section of the demand and supply side sample is given in Table \ref{tab:supply-demand}. Notice that it is a coincidence that the supply-side and demand-side reference subsamples are composed of exactly 104 respondents each. 
Of the 392 people interviewed, 24 belong to both the demand and supply samples, 80 belong only to the supply sample, 208 can only provide us with information useful for knowing the respondents' confidence towards the possible success of a crowd-shipping initiative.

\begin{table}
\centering
\begin{tabular}{rrrr}
							\toprule
               &       & \multicolumn{2}{r}{Living in Brescia} \\
               &       & \multicolumn{2}{r}{[Demand side]} \\[2mm]
               &       & {Yes} & {No} \\ 
							\midrule
 Metro user              & Yes   &  24 (6.12\%) & 80 (20.41\%) \\[2mm]
{[Supply side]} & No    &        80 (20.41\%) &   208 (53.06\%) \\
							\bottomrule
    \end{tabular}%
	\caption{Cross-section of demand and supply side of the sample.}
  \label{tab:supply-demand}%
\end{table}%

\subsubsection*{Demand side analysis}
To better understand the demand side of the sample, we collected information about the propensity to online shopping and the zip code of residence. The respondents living in Brescia municipality are 104, and are distributed by zip code as shown in Table \ref{tab:ZIP}.

\begin{figure}
  \centering
    \includegraphics[width=0.8\linewidth]{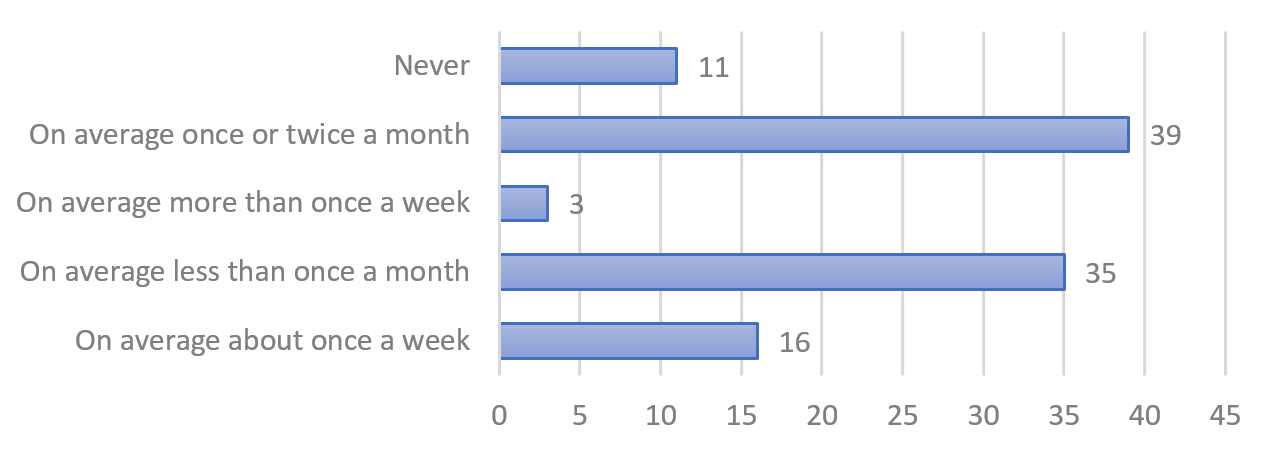}%
\caption{Frequency of online orders.}
\label{fig:demand-frequency}
\end{figure}

\begin{table}
  \centering
    \begin{tabular}{rrp{5mm}rr}
    {ZIP} & {\#}  & & {ZIP} & {\#} \\
		\toprule
    25123 & 27    & & 25122 & 4 \\
    25128 & 14    & & 25126 & 4 \\
    25121 & 13    & & 25135 & 3 \\
    25124 & 10    & & 25136 & 3 \\
    25125 & 8     & & 25134 & 2 \\
    25133 & 8     & & 25129 & 1 \\
    25127 & 6     & & 25131 & 1 \\
		\bottomrule
    \end{tabular}%
  \caption{Demand side: ZIP code distribution.}
  \label{tab:ZIP}%
\end{table}%
%
The distribution of online order frequencies in the sample of residents in Brescia is depicted in Figure \ref{fig:demand-frequency}. Notice that 56\% of the respondents order online at least once a month; 44\% of the respondents order online rarely or never. 
We posed the following question.
\begin{itemize}
\item \emph{Would you be willing to have the goods delivered to your home by an ``occasional shipper'', that is a person who does not do it for work but who has made himself available to make the delivery?} It turned out that among the 93 respondents who order online, 80 (86\%) would be willing to receive shipments from occasional shippers. The sample does not show relevant differences among age classes and gender. 
\end{itemize}
It would be interesting to analyze possible correlations between ZIP code, the inclination to online shopping, the willingness to be served by a crowd-shipper. However, this analysis would require larger numbers of respondents to be statistically sound.

\subsubsection*{Supply-side analysis}
\begin{figure}
  \centering
    \includegraphics[width=0.8\linewidth]{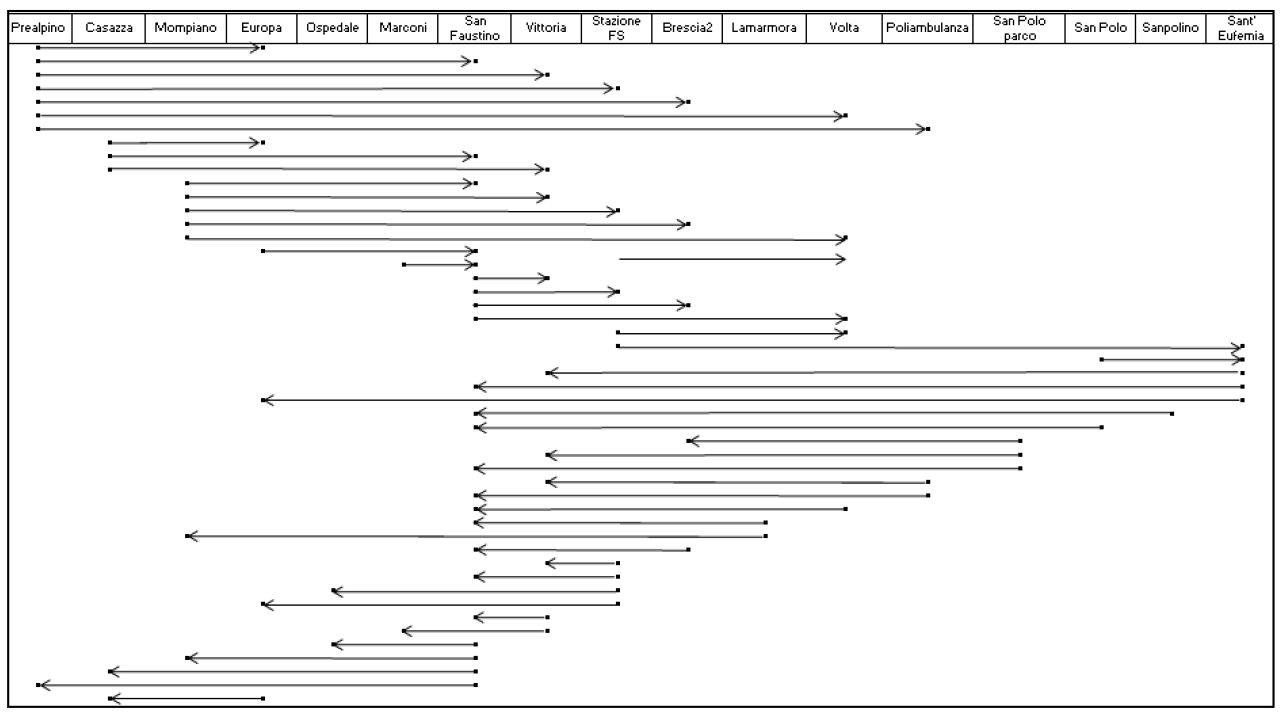}%
\caption{The first trips of the day in the sample.}
\label{fig:morning-trips}
\end{figure}

To analyze the supply side, we first ask to the 104 metro users how often they use the metro on average. The answers are reported in Table \ref{tab:frequenza-viaggi}. To visualize the variety of trips, in Figure \ref{fig:morning-trips} we represent the departure and arrival stations stated by the respondents for the first trip of the day (we assume the second trip is reversed). 
The typical trip described by frequent users of the metro involves departure from the stations of \emph{Prealpino}, \emph{Casazza} or \emph{Sant'Eufemia Buffalora} (located at the entrance to the city from different directions and close to large parking lots) or from \emph{Stazione FS} (the train station). The destination station of the typical trip is mainly located near \emph{Vittoria} or \emph{San Faustino}, two destinations in the city center, in areas characterized by car traffic limitations. 

\begin{table}
  \centering
    \begin{tabular}{rr}
    \toprule
    1 or 2 days per week & 21\% \\
    \midrule
    3 days per week & 21\% \\
    \midrule
    4 days per week & 22\% \\
    \midrule
    5 days per week & 27\% \\
    \midrule
    5 days per week & 9\% \\
    \bottomrule
    \end{tabular}%
	\caption{Frequencies of trips among metro users.}
  \label{tab:frequenza-viaggi}%
\end{table}%

We posed the following question to the metro users.
\begin{itemize}
\item \emph{During one or more of your metro rides, would you be willing to carry a package for free, thus becoming an occasional crowd-shipper? (Weight and size of a normal backpack).} It turned out that 69 respondents out of 104 (66\% of the total) would be willing to collaborate for free.
\end{itemize}
To the 35 respondents who were not available as crowd-shipper for free, we asked the following. 
\begin{itemize}
\item \emph{Would you be willing to carry the package for a reimbursement equal to recovering the costs of the trip?}
In this case, 19 out of 35 respondents (54\%) would be willing to cooperate as occasional transporters in exchange for reimbursement of the metro ticket.
\end{itemize}
To the 16 respondents not willing to collaborate in exchange for ticket reimbursement were finally offered compensation in excess of the cost of the ticket.
\begin{itemize}
\item \emph{Would you be willing to do it for more than the cost of the trip?} Here, 10 respondents (63\%) would accept to collaborate under these conditions. 6 respondents would not be willing under any circumstances.
\end{itemize}

To the 98 respondents who, possibly with compensation, say they would be willing to carry packages (the size of a backpack) we asked willingness to detour from the subway station. 
\begin{itemize}
\item \emph{If you were required to pick up/deliver the package outside of the Metro station, would you still accept?} The answers are summarized in Table \ref{tab:detour}.
\end{itemize}
\begin{table}
  \centering
    \begin{tabular}{p{18em}r}
    \toprule
    Yes, without extra charge, regardless of detour time & 2.0\% \\
    \midrule
    Yes, without extra charge, for a maximum 5 minute walking detour & 56.1\% \\
    \midrule
    Yes, without extra charge, for a maximum walking distance of 10 minutes & 13.3\% \\
    \midrule
    Yes, with an extra charge for the detour & 20.4\% \\
    \midrule
    No, in any case & 8.2\% \\
    \bottomrule
    \end{tabular}%
		 \caption{Accepted detour from the original trip of potential crowd-shippers.}
  \label{tab:detour}%
\end{table}%

\bigskip
We may summarize the analysis of the supply side as follows: about 66\% of the sample would be willing to collaborate free of charge; an additional 19\% of the sample would be willing to cooperate for a reimbursement equal to the cost of the ticket; an additional 9\% of the sample would be willing to cooperate for a reimbursement greater than the cost of the ticket; only 6\% of the sample would in no case be willing to collaborate in the initiative. 
In the light of these findings, our answer to research question Q1 is affirmative.

%

\section{Computational study}
\label{sec:prodes}
To assess the possible advantages of using crowd-shipping as an integration to conventional home delivery, we consider the following situation. Consider a store that manages home delivery of packages to customer located in Brescia. Every morning, a given set of customers has to be visited to deliver packages. Each delivery is carried out either by a van or by crowd-shipping. The packages assigned to the van are loaded on the vehicle at the store. The van leaves the store, visits all assigned customers to deliver the corresponding packages and, finally, goes back to the store. Each package assigned to crowd-shipping is picked up by the first available and compatible crowd-shipper, and delivered by him/her to the corresponding customer. A crowd-shipper is a metro passenger that picks up a package from the store (close to the departure metro station), travels by metro, and delivers the package to a customer's home located near the arrival metro station. Hence, only customers living sufficiently close to a metro station can be served by a crowd-shipper. To maintain the analysis simple, we assume that:
\begin{itemize}
\item the availability of crowd-shippers is such that any delivery assigned to crowd-shipping can find a corresponding shipper (unbounded availability of shippers);
\item there are no time restrictions for the deliveries (no time windows).
\end{itemize}
We assume that any delivery carried out by a crowd-shipper has the same fixed cost (the reward for the shipper), whereas the cost of the deliveries carried out by the van depend on the length of the route and the operating cost of the van.

In this section, our aim is to make a quantitative assessment of the advantages and disadvantages of crowd-shipping using metro passengers by optimizing the system described above under different scenarios. To do that, we need first to develop an optimization model, then to build a realistic case study based on the results of the questionnaire described above, and finally to run a computational study.

\subsection{A Vehicle Routing Problem with crowd-shippers}
\label{subsec:VRPOD}
In this section, we develop an optimization model to simulate the delivery system described above. 

Let $G = (N,A)$ denote a complete directed graph with node set $N = \{0, 1, \ldots, n\}$ and arc set $A = \{(i,j): i,j \in N, i \neq j\}$. Node 0 corresponds to the store, nodes $1, \ldots, n$ correspond to the daily set of customers to be served. Each arc $(i,j) \in A$ has a length $d_{ij}$ and a cost $c_{ij}$. Each customer $i \in N \setminus\{0\}$ has a given demand of $q_i$ units (measuring the size of the corresponding package). There are $m$ identical vehicles available, each one with a capacity of $Q$ units. 

Let $S \subseteq N \setminus \{0\}$ be the set of customers that can be served by a crowd-shipper. They are the customers located within a given maximum distance $\delta$ from a metro station, where $\delta$ is expressed for example as walking time. Let $V$ be set of customers that must be served by a standard vehicle, i.e., let $V = N \setminus (\{0\}\cup S)$.

A shipper receives a fixed reward $p$ as compensation for making a delivery to any customer $i \in S$. Vehicle routes need to satisfy a capacity constraint: the total demand of the customers served on a single route cannot exceed $Q$. It is implicitly assumed that a shipper can accommodate the demand $q_i$ of customer $i \in S$ he/she is willing to serve. The cost of a vehicle route $r$ is the  sum of the costs of the arcs in the route, i.e., $\sum_{(i,j) \in r} c_{ij}$. From the company's perspective, the objective is to minimize the total costs, that is the sum of the costs incurred by the regular vehicles and the cost incurred for compensating the crowd-shippers. 
Note that if there are no crowd-shippers, the problem just described becomes the standard Capacitated Vehicle Routing Problem (CVRP). 

Let $x_{ij}$ be a binary variable indicating whether a vehicle traverses arc $(i, j)$. Let $y_{ij}$ indicate the load a vehicle carries on the arc $(i, j)$. Let $z_i$ be a binary variable indicating whether customer $i \in S$ is visited by a vehicle or a crowd-shipper. More precisely, we assume $z_i = 1$ if $i$ is visited by a vehicle; $z_i = 0$ if $i$ is visited by a crowd-shipper. We formulate the following \emph{Vehicle Routing Problem with Crowd-Shippers} (VRPCS):

\begin{alignat}{3}
     \min &  \sum_{(i,j) \in A} c_{ij} x_{ij} + p \sum_{i \in C} (1-z_i) \label{VRPS-obj} \\
\subjto{} & \sum_{j| (i,j) \in A} x_{ij} = \sum_{j| (j,i) \in A} x_{ji} = 
            \begin{cases}
						   1   & i \in V \\
							 z_i & i \in S 
						\end{cases} \label{VRPS-assignment}\\
          & \sum_{j| (0,j) \in A} x_{0j} = \sum_{j| (j,0) \in A} x_{j0} \leq m \label{VRPS-routenum}\\
          & \sum_{j| (j,i) \in A} y_{ji} - \sum_{j| (i,j) \in A} y_{ij} =
					\begin{cases}
						   q_i     & i \in V \\
							 q_i z_i & i \in S
						\end{cases} \label{VRPS-vehicleflow1}\\
          & \sum_{j| (j,0) \in A} y_{j0} - \sum_{j| (0,j) \in A} y_{0j} = 
						-\left(\sum_{i \in V} q_i + \sum_{i \in S} q_i z_i\right) \label{VRPS-vehicleflow0} \\
          &  y_{ij} \leq Q x_{ij}  & (i,j) \in A \label{VRPS-capacity}\\
					&  y_{i0} = 0 & i \in N \setminus \{0\} \label{VRPS-emptyvehicle} \\
          &  x_{ij} \in \{0,1\} & (i,j) \in A \label{VRPS-xbinary} \\
          &  z_{i} \in \{0,1\} & i \in S \label{VRPS-zbinary} \\
					&  y_{ij} \geq 0 & (i,j) \in A \label{SSCFL-ynonneg}
\end{alignat}
The objective function (\ref{VRPS-obj}) aims at minimizing the total cost. Constraints (\ref{VRPS-assignment}) and (\ref{VRPS-routenum}) are flow conservation constraints. In particular, constraint (\ref{VRPS-routenum}) compels that no more than $m$ vehicles can be used. Constraints (\ref{VRPS-vehicleflow1}) and (\ref{VRPS-vehicleflow0}) ensure that demand of the customers visited by vehicles is satisfied and that subtours  are  prevented. Constraints (\ref{VRPS-capacity}) ensure that vehicle capacity is respected. Finally, constraints (\ref{VRPS-emptyvehicle}) establish that the vehicles must return empty to the depot. 

Though developed from a different setting, the above model is similar to the VRP with Occasional Drivers (VRPOD) \cite{2016-ArcSavSpe}.

\subsection{Brescia case study}
\label{subsec:Brescia}
In this section, evaluate empirically the potentials of crowd-shipping based on metro commuters in the city of Brescia, Italy. 
We first describe the case and the experimental setting, then the results and their discussion.

\subsubsection{Case and experimental setting}
\begin{figure}%
\begin{center}
\includegraphics[width=0.7\columnwidth]{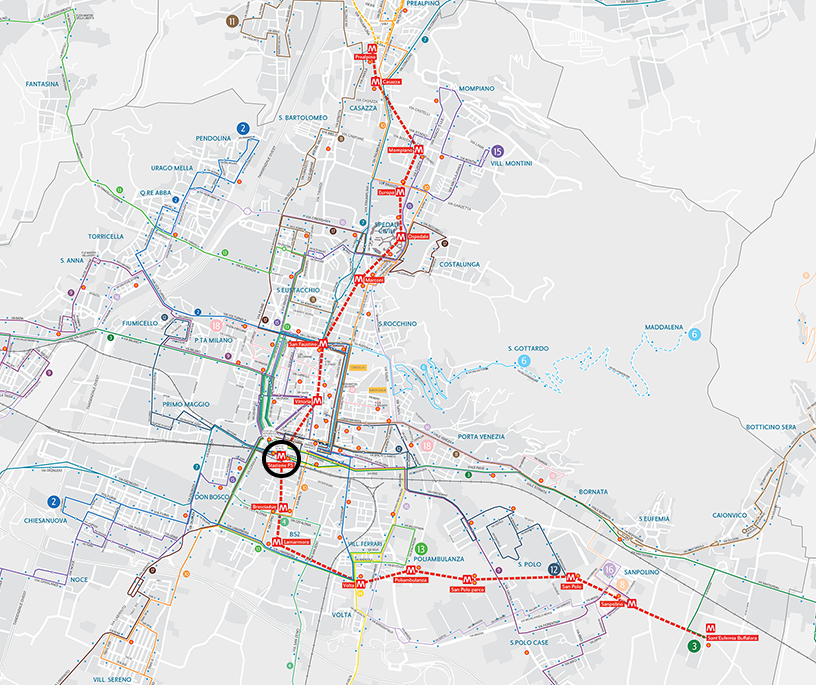}%
\end{center}
\caption{Map of Brescia metro. The line and the stations are represented in red. The main station is circled in black.}%
\label{fig:MetroMap}%
\end{figure}

We consider a delivery company that delivers packages in Brescia municipality starting from a depot located near the railway station, where also the main metro station is located. The Brescia metro is composed by a single line with an L-shape, with a terminus in the south-east part of the town and the other terminus in the northern part. The main station is adjacent to the railway station, and is positioned in the middle of the metro line (see Figure~\ref{fig:MetroMap}).

\begin{figure}%
\begin{center}
\includegraphics[width=0.7\columnwidth]{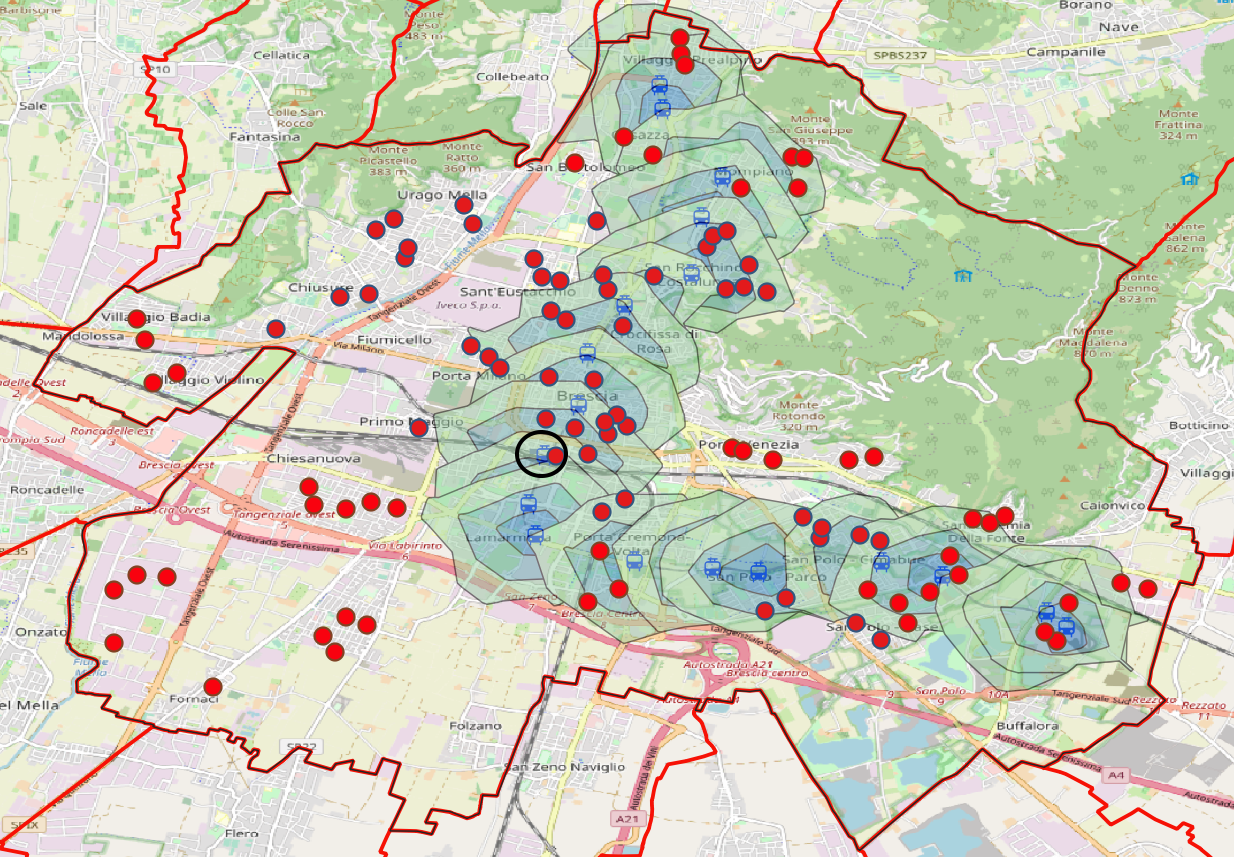}%
\end{center}
\caption{Customers (red circles) and isochrones at 5, 10, 15 minutes of walking distance from metro stations. The circled position corresponds to the main metro station and the depot.}%
\label{fig:CustIsoMap}%
\end{figure}

We consider 100 customers to be served within Brescia municipality. Given an integer $\delta$, we assume that a customer located at no more that $\delta$ minutes of walking distance from a metro station may be visited either by a vehicle or by a crowd-shipper. A customer located at more that $\delta$ minutes of walking distance from any metro station must be visited by a vehicle. Notice that the $\delta$ minutes can be interpreted as minutes walked by the crowd-shipper to reach a customer and deliver a package or as minutes walked by the customer to reach a locker placed in a metro station where a crowd-shipper has left the package. This distinction is irrelevant for the optimization model discussed above. The location of the customers and sample isochrones around the metro stations are illustrated in Figure~\ref{fig:CustIsoMap}.
The location of the customers is chosen to have a spacial distribution roughly correlated to the population density.

Given this setting, we test the following grid of parameters:
\begin{itemize}
\item \emph{Demand}: for all $i \in V \cup S$, $q_i$ is an integer randomly chosen in $\{1, 2, 3, 4, 5\}$; this should account for the differences in size of the packages. To consider different situations, we generate independently five different vectors of elements $q_i$, $i \in V \cup S$.
\item \emph{Vehicle capacity}: we consider $Q \in \{360, 180, 90, 45\}$. Given the demand distribution, these capacities should identify situations where the expected number of vehicles to be used is respectively 1, 2, 4, 8. In any case, we set $m = 10$. 
\item \emph{Cost coefficients and price of crowd-shipping}: we take as unit cost the cost per kilometer of a vehicle. Thus, for all $(i,j) \in A$, $c_{ij}$ is the road distance from $i$ to $j$ in kilometers. The price of crowd-shipping is expressed in the same unit, and we test four different prices, with $p \in \{0, 0.5, 1, 1.5, 2, 5\}$. For $p = 0$ (crowd-shippers offer their service for free), crowd-shipping should always be the choice when possible; for $p = 5$ (each crowd-shipper requires a reward corresponding to a 5-km run of a vehicle), crowd-shipping is always too expensive.
\item \emph{Customers compatible with crowd-shipping}: we assume that only the customers close to a metro station are compatible with crowd-shipping. We test three levels of acceptable walking distance $\delta$:  5, 10, and 15 minutes.
\end{itemize}
The total number of runs is thus $5 \times 4 \times 6 \times 3 = 360$.

We derived all road distances and the isochrones of walking distance from the metro stations using the OR tools in QGIS3; the optimization model was implemented in Python 3.8 using Gurobi 9.1 as solver. Since the purpose of this analysis is to get insights on the managerial problem and not to design the efficient solution of model VRPCS, we set for Gurobi a time limit (parameter \texttt{TimeLimit}) of 600 seconds and a MIP gap tolerance (parameter \texttt{MIPGap}) of 1\% on each instance. All the other Gurobi parameters are set to default. The tests have been performed on a 64-bit Windows laptop with Intel Core processor i5-8250U, 1.8GHz, with 8 GB of RAM.


\subsubsection{Results and Discussion}
We first evaluate the computational effort required by model VRPCS and the quality of the obtained solution, then we discuss the managerial insights obtained from our simulation. In both cases, we resort mostly on a graphical representation of the results. 

\begin{figure}%
\begin{center}
\includegraphics[width=0.7\columnwidth]{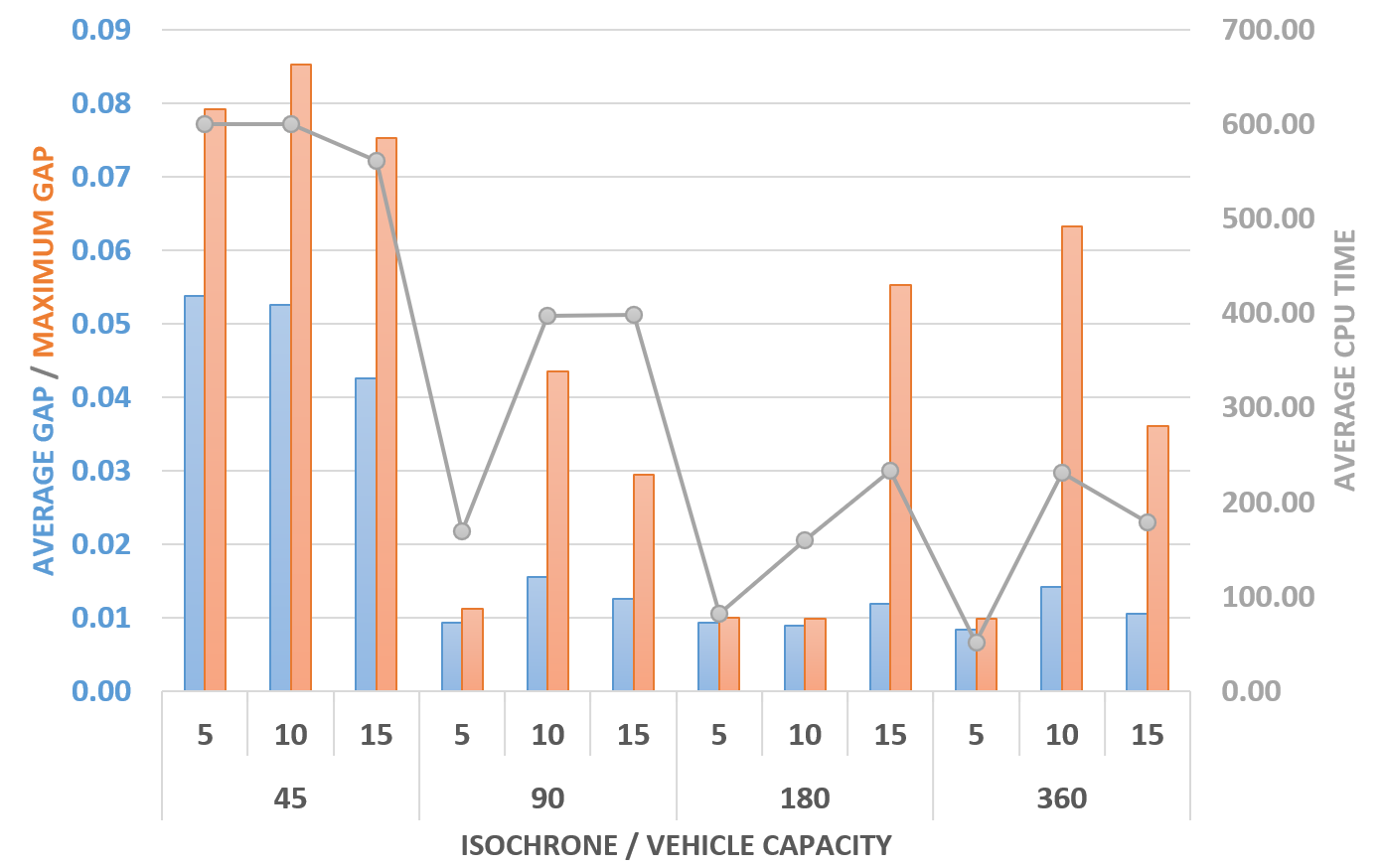}%
\end{center}
\caption{Optimality gap and CPU time as functions of crowd-shipping compatibility and vehicle capacity.}%
\label{fig:Q+Iso-vs-gap+time}%
\end{figure}

Concerning the computational effort, in Figure \ref{fig:Q+Iso-vs-gap+time} we plot the average and maximum optimality gap and the average CPU time as functions of the vehicle capacity $Q$ and the isochrone (i.e., the number of customers compatible with crowd-shipping). The chart shows that the hardness of the problem decreases with the capacity of the vehicle, i.e., it is inversely proportional to the average number of vehicles needed to serve the customers. In particular, while there is an almost regular decreasing trend in the CPU time, the decrease in gap appears to be step-wise, with a large improvement from $Q=45$ to $Q=90$ followed by a quite stable behavior in the following capacities. This could be explained by the fact that the time limit of 600 seconds implied that most of the instances with $Q=45$ have been terminated prematurely. On the other hand, the role of the number of customers compatible with crowd-shipping is not apparent.

\begin{figure}%
\begin{center}
\includegraphics[width=0.7\columnwidth]{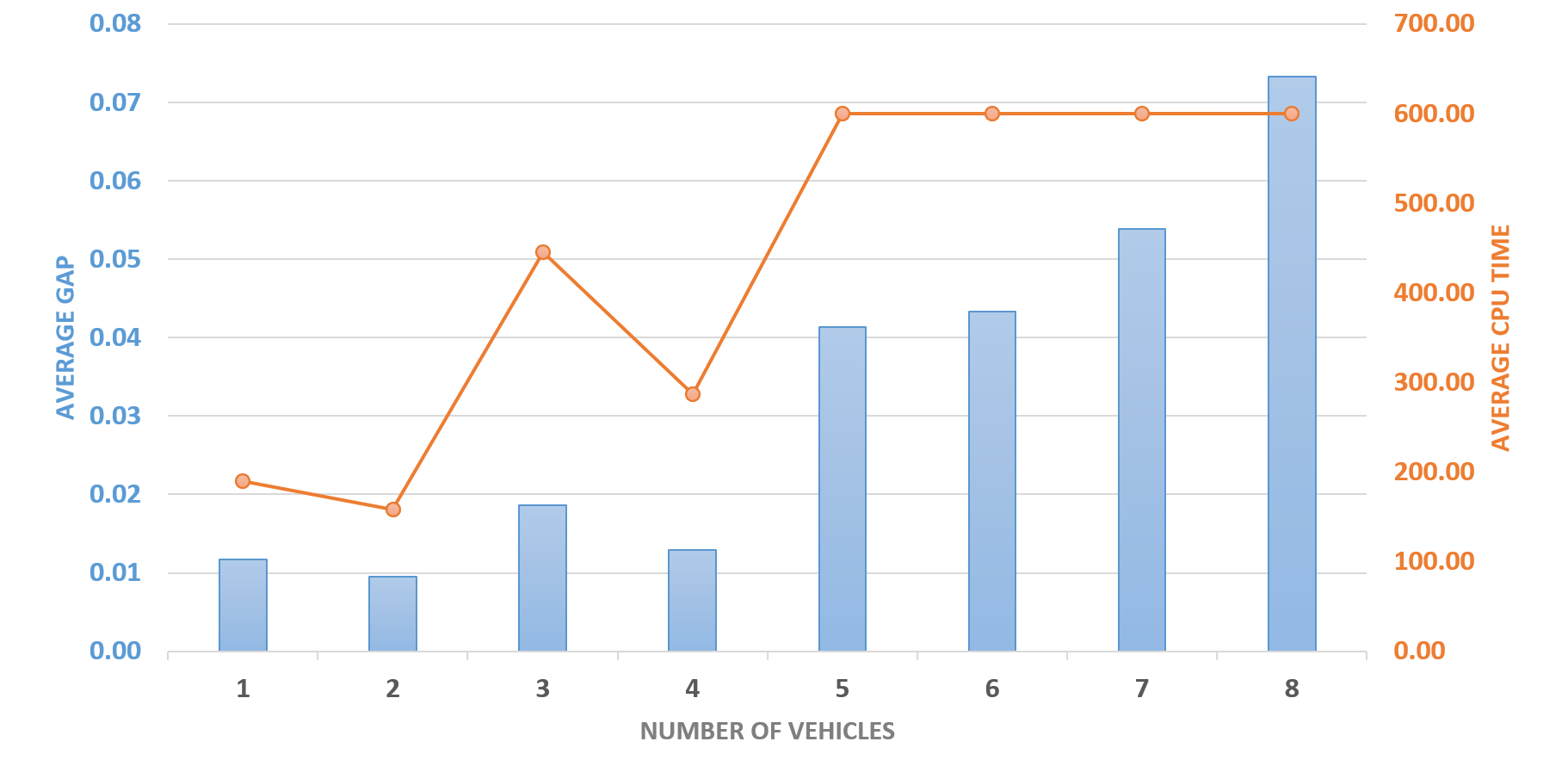}%
\end{center}
\caption{Optimality gap and CPU time as functions of crowd-shipping reward.}%
\label{fig:vehicle_num-vs-gap+time}%
\end{figure}

The strict correlation between number of used vehicles and hardness of the problem is confirmed by Figure \ref{fig:vehicle_num-vs-gap+time}. We see that, with five or more vehicles, all runs are terminated by the time limit. 

\begin{figure}%
\begin{center}
\includegraphics[width=0.8\columnwidth]{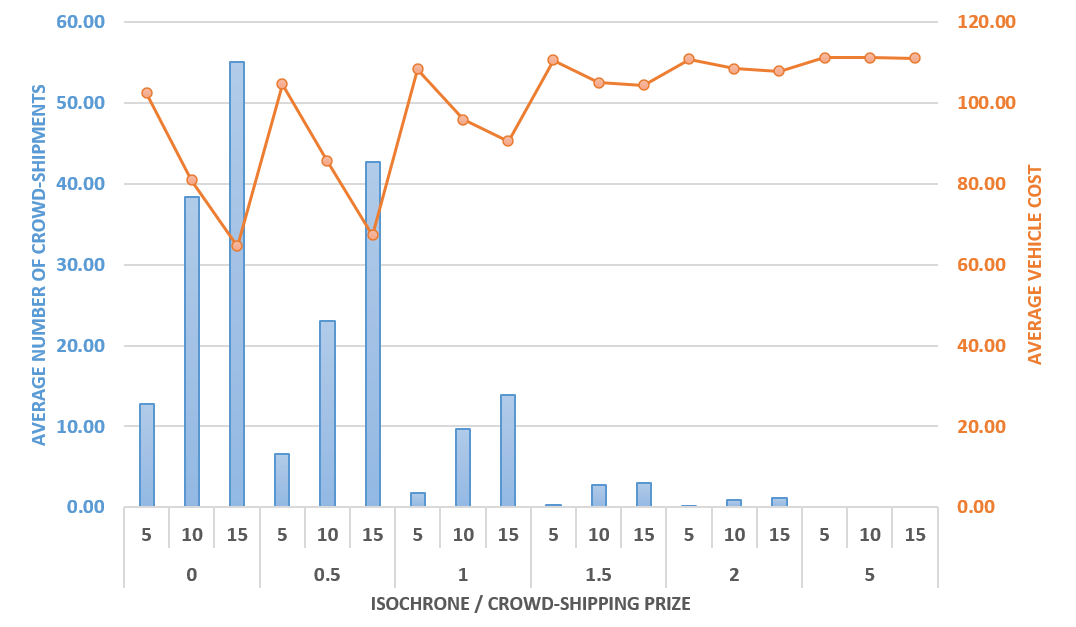}%
\end{center}
\caption{Number of crowd-shipments and vehicle cost as functions of potential shipmnets and crowd-shipping reward.}%
\label{fig:p+iso-vs-crowd_num+vehicle_cost}%
\end{figure}

\begin{figure}%
\begin{center}
\includegraphics[width=0.7\columnwidth]{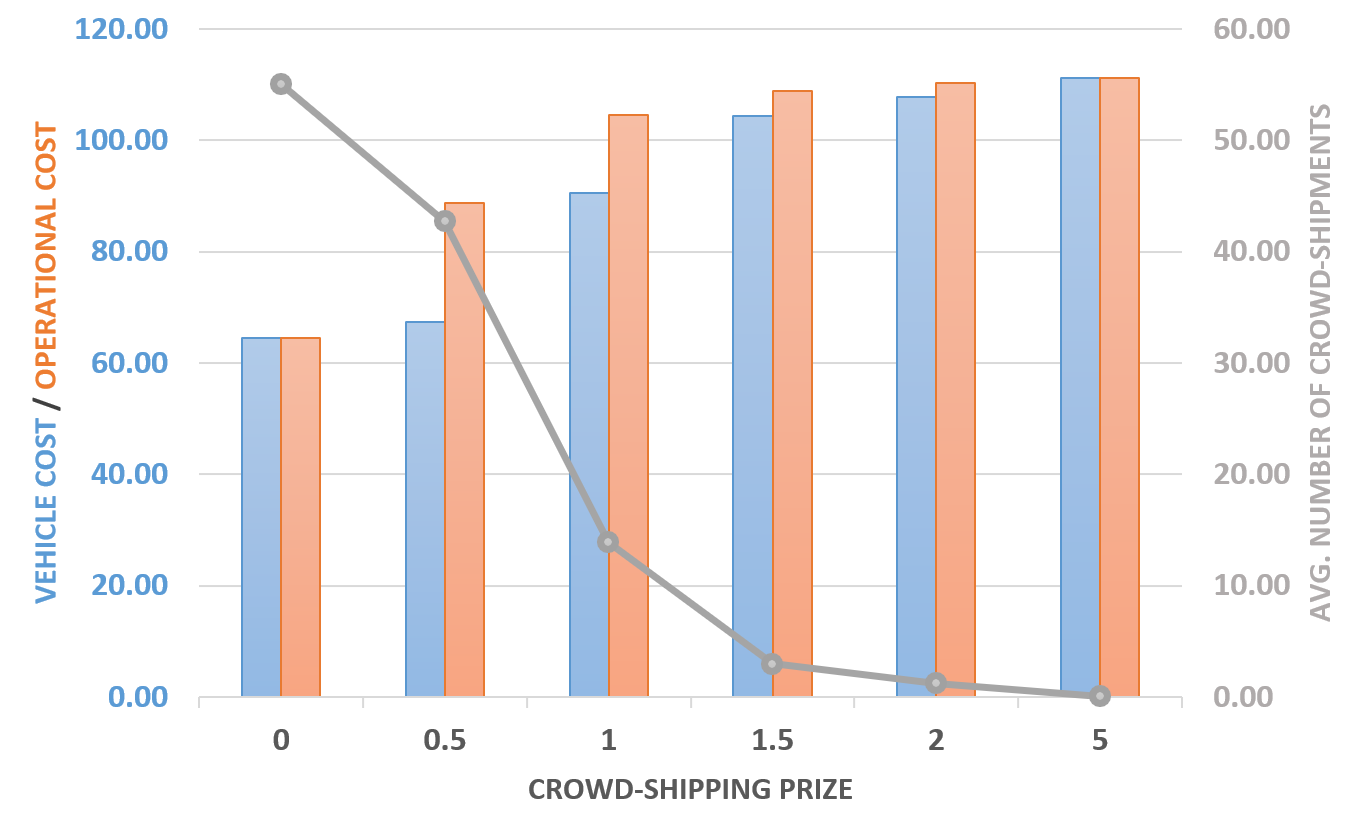}%
\end{center}
\caption{Costs and number of crowd-shipments as functions of crowd-shipping reward.}%
\label{fig:prize-cost&crowd}%
\end{figure}

Turning to the insights that model VRPCS gives in our case study, we analyze in Figure \ref{fig:p+iso-vs-crowd_num+vehicle_cost} how the number of potential crowd-shipments and the crowd-shipping reward impact the average number of achieved crowd-shipments and the average traveling cost of the vehicles. As one may expect, there is a strong negative correlation between the number of customers compatible with crowd-shipping and the traveling cost of the vehicles: the higher the chance to serve customers by crowd-shipping, the lower the usage of vehicles. Moreover, there is a positive correlation between the crowd-shipping reward and the traveling cost for vehicles: the higher the cost for crowd-shipping, the higher the distance traveled by vehicles. Moreover, the reward paid to a crowd-shipper has a huge impact on the convenience of resorting to crowd-shipping. The latter relation is better analyzed in Figure \ref{fig:prize-cost&crowd}, where we plot the average behavior of the solutions when the price of each crowd-shipment is increasing. As we may expect, the higher the price, the lower the convenience to use crowd-shipping and, beyond a certain price, crowd-shipping is not convenient at all. Concerning the operational costs (traveling of vehicles plus crowd-shipping), the average total expense when $p=5$ (no crowd-shipping used) is 111.18, whereas the the average total expense when $p=0$ (maximum usage of crowd-shipping) is 82.67, with an average saving of 25.64\%. For opposite reasons, when $p=5$ and when $p=0$, the operational cost is entirely due to standard vehicles. In Table \ref{tab:pippo}, we analyze the trade-off between vehicle cost and crowd-shipping cost for the entire range of rewards for crowd-shipping. For high reward values, the cost saving given by crowd-shipping is negligible. However, when the compensation for crowd-shipping is low, the saving is sensible. In particular, when $p=0.5$ (i.e., when the compensation is equal to the cost of a deviation of 0.5 km for a vehicle) we have an average total saving of almost 12\%, with a saving in traveling cost of 22.76\%. The last columns also gives an interesting insight: when $p=0.5$, the average number of required vehicles is 2.70 whereas, without crowd-shipping, the average number of required vehicles is 3.63. That means that, on average, on the Brescia case we could save one vehicle out of four, with an additional 25\% saving in amortization cost.


\begin{table}
  \centering
    \begin{tabular}{rrrr}
		\toprule
    \textbf{$p$} & \multicolumn{1}{p{7.0em}}{Avg. Traveling Cost Saving} & \multicolumn{1}{p{8.0em}}{Avg. Operational Cost Saving} & \multicolumn{1}{p{7.5em}}{Avg. Number of Vehicles} \\
    \midrule
    0.0   & 25.64\% & 25.64\% & 2.57 \\
    0.5   & 22.76\% & 11.91\% & 2.70 \\
    1.0   & 11.57\% & 3.96\% & 3.17 \\
    1.5   & 4.05\% & 1.39\% & 3.45 \\
    2.0   & 1.87\% & 0.55\% & 3.57 \\
    5.0   & 0.00\% & 0.00\% & 3.63 \\
    \bottomrule
    \end{tabular}%
		  \caption{An analysis of costs and number of required vehicles with respect to crowd-shipping cost.}
  \label{tab:pippo}%
\end{table}%

\begin{figure}
  \centering
  \begin{subfigure}{0.475\textwidth}
    \centering
    \includegraphics[width=0.8\linewidth]{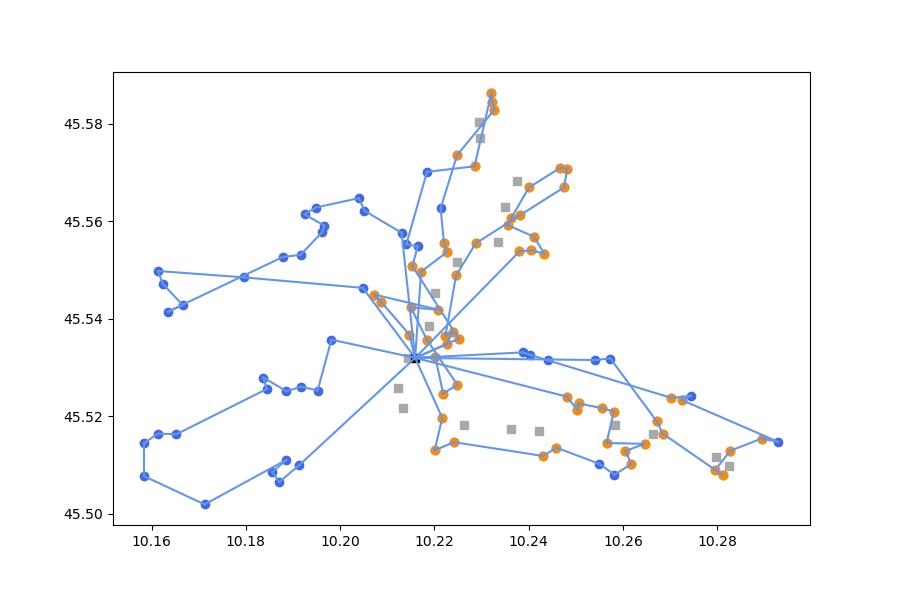}%
		\caption{$p = 5.0$ \label{subfig:f}}
  \end{subfigure}
  \hfill
  \begin{subfigure}{0.475\textwidth}
    \centering
    \includegraphics[width=0.8\linewidth]{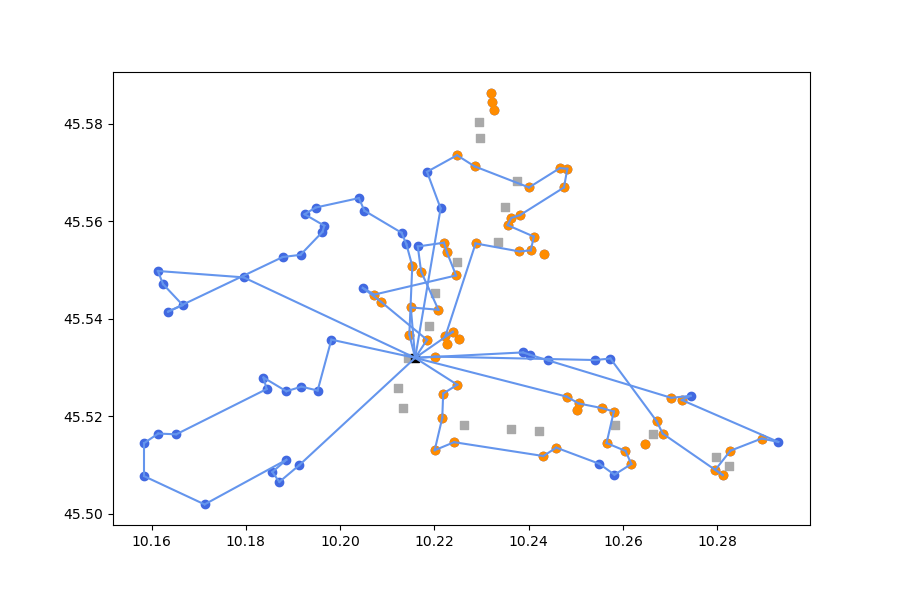}%
		\caption{$p = 1.5$ \label{subfig:d}}
  \end{subfigure}
	\\
	\begin{subfigure}{0.475\textwidth}
    \centering
    \includegraphics[width=0.8\linewidth]{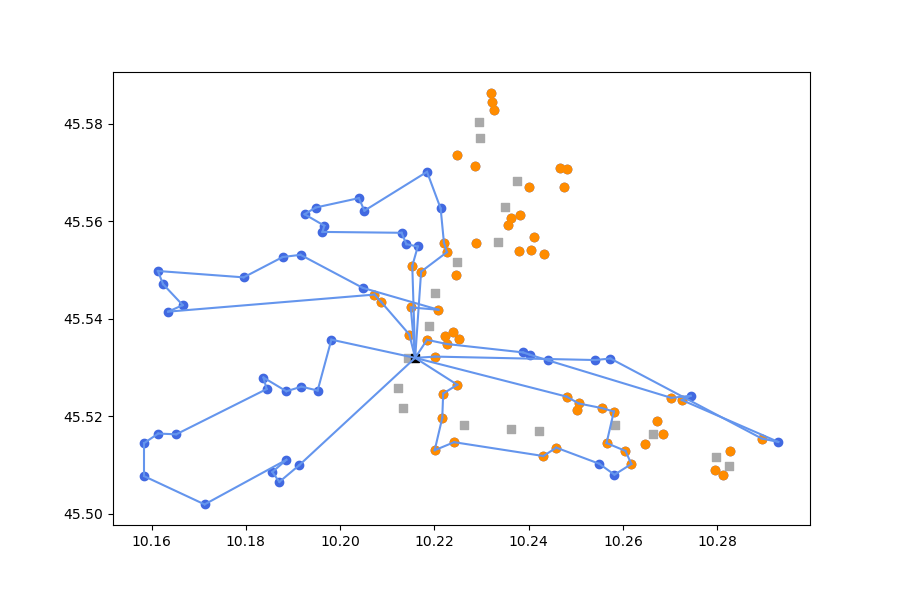}%
		\caption{$p = 1.0$ \label{subfig:c}}
  \end{subfigure}
  \hfill
  \begin{subfigure}{0.475\textwidth}
    \centering
    \includegraphics[width=0.8\linewidth]{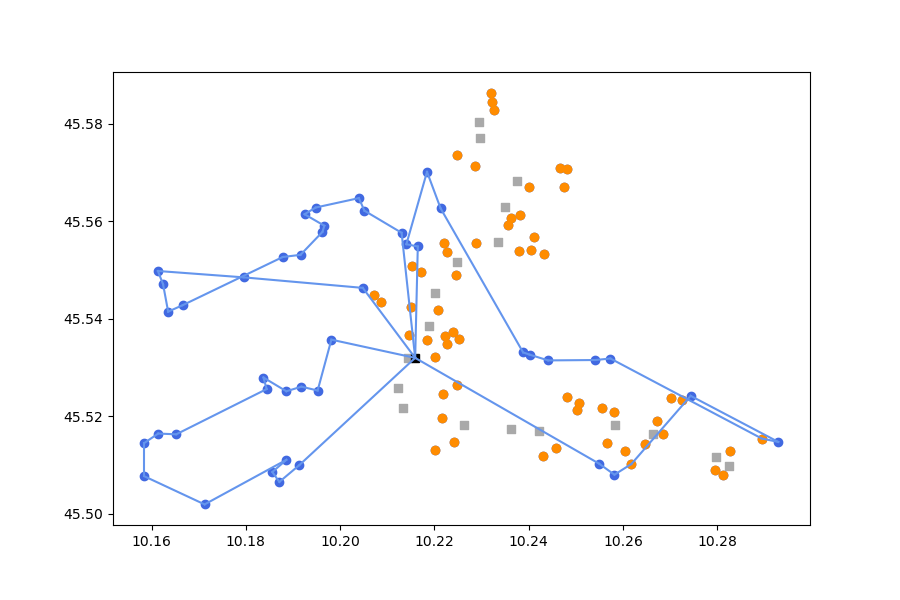}%
		\caption{$p = 0.0$ \label{subfig:a}}
  \end{subfigure}
\caption{Evolution of the routes when the reward for crowd-shipping decrease.}
\label{fig:dummy}
\end{figure}

Finally, in Figure \ref{fig:dummy} we give a graphical view of the routes when crowd-shipping becomes cheaper. More precisely, we plot the solutions obtained under the first demand scenario when $Q=45$ and 15 minutes of walking distance is allowed for crowd-shipping. We consider different decreasing values of crowd-shipping reward: $p=5.0$ (no crowd-shipping), $p=1.5$, $p=1.0$, $p=0.0$ (free crowd-shipping). In Figure \ref{fig:dummy}, the depot is the black square, the metro stations are the gray squares, the customers compatible with crowd-shipping are the orange dots, and the customers not compatible are the blue dots. In this particular instances, crowd-shipping could dramatically reduce the number of needed routes, passing from 7 to 3. Though a large part of the customers cannot be reached by crowd-shipping and require a standard delivery service, it is apparent that the availability of cheap crowd-shipping may greatly simplify the delivery work for a package carrier.

In the light of the above findings, our answer to research question Q2 is also affirmative.

\section{Conclusions}
\label{sec:concl}
We have analyzed the potentialities of crowd-shipping based on public transport for last-mile delivery in a medium-size European city. In particular, we focused in Brescia city, in Northern Italy. Our study is twofold: on one side, through a survey, we have investigated the interest of metro users to act as crowd-shippers; on the other side, we have developed a computational study to simulate the opportunities for a shipper with a depot located near the main metro station and the possibility to serve part of its customers by using crowd-shippers. Essentially, the results are that: (a) there is a common awareness that sustainable delivery modes are crucial to improve the quality of life in our cities and a substantial willingness to cooperate, even with no monetary reward; (b) for a carrier, the implementation of a crowd-shipping initiative may have not only a useful impact in terms of image, but also a remarkable impact in terms of operational costs.

It should be noted that our study shows positive results even if we consider a case study where just one line of public transport can be exploited and hence only a limited portion of shipments can be fulfilled by crowd-shipping. In cities where a network composed by different lines can be used, we would expect even more positive results, though the coordination of crowd-shippers may become more complex to manage. 

The optimization model we use assumes that there is availability of crowd-shippers for any shipment to a compatible delivery location. This assumption is reasonable if the number of shipments is very low with respect to the number of public transport commuters. 
In this case, we may assume there is always a commuter starting from the pickup point (the central station in our case study) and ending at the destination point. However, for large cities with a complex public transport system, this assumption may become questionable, and the model may require an extension, for instance considering a limited availability of crowd-shippers for every  origin-destination pair.

Finally, on the algorithmic side, a regular use of the optimization model would require the development of an efficient heuristic to get good solutions quickly and to revise plans in real time.

\subsection*{Declarations}

\noindent
\textbf{Conflicts of interest} The authors declare that they have no conflict of interest.

\medskip
\noindent
\textbf{Availability of data} Upon request to the corresponding author.

\medskip
\noindent
\textbf{Code availability} Upon request to the corresponding author.

\bibliography{crowd}

\end{document}